\documentstyle[12pt]{article}
\newcommand{\beq}{\begin{equation}}
\newcommand{\eeq}{\end{equation}}
\newtheorem{theorem}{Theorem}
\newtheorem{proposition}[theorem]{Proposition}

\begin{document}

\centerline{\large\bf CHAOS IN A CONTINUOUS-TIME}
\centerline{\large\bf BOOLEAN NETWORK}
\bigskip
\centerline{R. Edwards}
\medskip
\centerline{Department of Mathematics and Statistics }
\centerline{University of Victoria }
\centerline{P.O.Box 3045, STN CSC, Victoria, BC, Canada V8W 3P4 }
\bigskip\smallskip
\noindent{\bf Abstract:} 
Continuous-time systems with switch-like behaviour occur in chemical kinetics, gene regulatory networks and neural networks. Networks with hard switching, as a limiting case of smooth sigmoidal switching, retain the richest possible range of behaviors but are mathematically more tractable. The form of an underlying discrete (fractional-linear) map encodes information on existence, stability and exact periods of periodic orbits. In richly connected structures with four or more variables, aperiodic behaviour can occur. We investigate a simple 4-dimensional example with Boolean interaction terms in which a Smale horseshoe-like object reveals chaotic dynamics.

\smallskip
\noindent{\bf AMS(MOS) subject classification:} 34C35, 92B20, 94C10

\bigskip
\centerline{\bf 1. INTRODUCTION}
\medskip
\noindent Aperiodic behavior in systems of three or more ordinary differential equations (ODEs) often requires careful tuning of parameters and there are few general principles that enable us to say when a given system can have such behaviour. Here we investigate a class of systems of ODEs in which aperiodic behavior results from a particular structure of coupling between the variables rather than finely tuned parameter values. 

For systems of interacting quantities that vary continuously but are dominated by switch-like behavior, Glass [\ref{g75a},\ref{g75b},\ref{g77}] proposed an approach, using structural equivalence classes based on state transition diagrams on an $n$-cube, to aid in relating dynamical behaviour to the structure of interactions in the system. Glass \& Pasternack [\ref{gp78a},\ref{gp78b}] then showed that periodic oscillations of various types occur with particular structural classes. This approach is exact in the case of `hard' switching (interaction terms depend only on whether the other variables are above or below threshold), but appears to be a good model for the case of steep sigmoidal switching. The methods of analysis work only for the case of identical decay rates for all variables.

These networks were first investigated (Glass [\ref{g75a},\ref{g75b},\ref{g77}]) in the context of chemical and biological oscillations: kinetics of interacting chemical species, interacting biological species, gene and enzyme regulatory networks and neural networks. We have recently argued that networks of this type may underlie transitions between irregular and regular tremor generated by the brain's motor circuitry in Parkinson's disease (Edwards {\em et al} [\ref{ebg99}]). As a consistent body of theory has emerged in the case of identical decay rates and hard switching, and as the existing nomenclature is confusing, we propose to call these systems `Glass networks' after their originator. The hard switching makes them remarkably tractable, yet they remain very rich in dynamical possibilities. 

The structural equivalence classes have a particularly simple representative in which the values of the interaction terms are just Boolean functions (`Boolean Glass networks'). It is not clear to what extent these are representative of the dynamics of their respective classes, but they do permit complex behavior.

Aperiodic behavior appears to be common in Glass networks with many (6 or more) variables (Lewis \& Glass [\ref{lg92}]; Mestl {\em et al} [\ref{mbg97}]; Glass \& Hill [\ref{gh98}]; Edwards {\em et al} [\ref{ebg99}]) and chaos was found in one Glass network of 4 variables with a particular set of parameters (Mestl {\em et al} [\ref{mlg96}]). However, numerical simulations suggest that aperiodic behavior is not so rare even among 4-dimensional Boolean Glass networks, and here we show how available techniques can be used to analyze the behavior of one such network, via an object resembling the Smale horseshoe.

\bigskip
\centerline{\bf 2. GLASS NETWORKS} 
\medskip
\noindent A Glass network is a system of the form
\beq\dot{y_i}=-y_i + F_i (\tilde{y}_1,\tilde{y}_2,\ldots ,\tilde{y}_n)\,,\quad i=1,\ldots ,
n\,,\label{eq:network}\eeq
where 
\beq\tilde{y}_i=\left\{\begin{array}{ll}
0 & \mbox{if $y_i<0$} \\
1 & \mbox{if $y_i>0$}
\end{array} \right.\,.\label{eq:step}\eeq
Note that all thresholds are 0 (Equation~\ref{eq:step}), and $F_i$ depends only on the signs of the variables $y_i$. Systems with non-zero thresholds and values of $\tilde{y}_i$ other than 0 and 1, possibly different for each variable, can be reduced to the above form by appropriate transformations. Similarly, a decay rate parameter may be put in the equations, and this may even depend on the $\tilde{y}_i$, but the analysis below will only apply if they are uniform among the variables of the system. This network can be considered a limiting case of smooth networks in which the step function (Equation~\ref{eq:step}) is replaced by a sigmoid, $\tilde{y}_i=g(y_i)$, 
whose range is the open interval $(0,1)$. 

Since for a Glass network there are a finite number of values $F_i$ ($n2^n$ of them), it is clear that solutions are globally bounded. In the sigmoidal case, we would need the additional assumption that $F_i$ is bounded on $(0,1)^n, \forall i$. Boolean Glass networks are defined as Glass networks for which $F_i(\tilde{\bf y})=\pm 1$ for all $i$ and all $\tilde{\bf y}$, so that the interactions are Boolean functions (note that we could equivalently have made the input values, $F_i\in \{0,1\}$ if the thresholds had been set at $\frac{1}{2}$). All the theory we will use applies to the general Glass networks, but our chaotic example will Boolean.

While weaker conditions suffice for some of the following theory, we will assume:
\begin{description}
\item[\bf Condition 1:] $F_i\ne 0\,,\quad \forall i\,,\forall \tilde{\bf y}\,,\quad$ and
\item[\bf Condition 2:] $F_i(\tilde{y}_1,\ldots \tilde{y}_i=0,\ldots ,\tilde{y}_n)=F_i(\tilde{y}_1,\ldots \tilde{y}_i=1,\ldots ,\tilde{y}_n)$,
\end{description}
where we use $\tilde{\bf y}$ for the vector $(\tilde{y}_1,\tilde{y}_2,\ldots ,\tilde{y}_n)^{\prime}$, and similarly use bold face for other vectors, the $^{\prime}$ denoting matrix transposition. Condition 2 states that $F_i$ does not depend on $\tilde{y}_i$, {\em i.e.}, that there is no self-input in the network.

\bigskip
\centerline{\bf 3. CYCLES AND PERIODIC ORBITS} 
\medskip
\noindent The analysis of Glass networks was begun by Glass \& Pasternack [\ref{gp78b}] and was further developed mainly by Mestl {\em et al} [\ref{mpo95}] and Mestl {\em et al} [\ref{mlg96}]. What follows is a brief summary of this work though some results are new. 

The main property of Glass networks that makes them tractable is that trajectories are piecewise-linear. For ${\bf y}=(y_1,y_2,\ldots ,y_n)$ in one orthant of phase space (and therefore with one fixed sign structure) the solution to Equation~\ref{eq:network} (in vector form) is
\beq{\bf y}(t)={\bf f}+({\bf y}(0)-{\bf f})e^{-t}\,,\label{eq:solnpiece}\eeq
which describes exponential approach to ${\bf f}=(f_1,f_2,\ldots ,f_n)={\bf F}(\tilde{\bf y})$ in a straight line. Thus, each orthant of phase space (with sign structure $\tilde{\bf y}$) has an associated focal point, ${\bf f}$, somewhere in ${\bf R}^n$. If trajectories in an orthant are directed to a focal point ${\bf f}$ within that orthant then once the orthant is entered no further switchings take place and ${\bf f}$ is a stable fixed point of the network dynamics. Otherwise, trajectories are formed of piecewise-linear segments between orthant boundaries, with sharp corners at the boundaries. Under Condition 2, there is no ambiguity in the direction of flow across an orthant boundary so trajectories are well defined there.

We now denote by ${\bf y}^{(k)}$ the $k^{th}$ such orthant boundary crossing on a trajectory and assume that ${\bf f}^{(k)}$, the focal point associated with the orthant being entered, does not lie in that orthant. The map from one boundary to the next can be represented as an operator ($M^{(k)}:{\bf R}^n\rightarrow {\bf R}^n$):
\beq{\bf y}^{(k+1)}=M^{(k)}{\bf y}^{(k)}=\frac{B^{(k)}{\bf y}^{(k)}}{1+\langle \psi^{(k)},{\bf y}^{(k)}\rangle } \,,\quad B^{(k)}=I-\frac{{\bf f}^{(k)}{\bf e}_j^{\prime}}{f_j^{(k)}}\,,\quad \psi^{(k)}=\frac{-{\bf e}_j}{f_j^{(k)}}\,,\label{eq:flm1}\eeq
where $j$ is the variable that switches at the $k^{th}$ step, ${\bf e}_j$ denotes the standard basis vector in ${\bf R}^n$ and the angle brackets denote the Euclidean inner product ($\langle\psi ,{\bf y}\rangle=\psi^{\prime}{\bf y}$).
Thus, $M^{(k)}$ is a fractional-linear map with a vector numerator and scalar denominator. The composition of such maps is again a fractional-linear map of the same form. Also, since these maps are between orthant boundaries where one of the $y_i$'s is always 0, they can be reduced by one dimension, by removing the appropriate row and column in each $B^{(k)}$, ${\bf y}^{(k)}$ and $\psi^{(k)}$. For a cycle, (a trajectory that returns to its initial orthant boundary), we arrive at (dropping the superscripts)
\beq M{\bf y}=\frac{A{\bf y}}{1+\langle\phi ,{\bf y}\rangle }\,,\label{eq:Acycle}\eeq
where $A$ is $(n-1)\times (n-1)$, $\phi\in {\bf R}^{n-1}$ and ${\bf y}\in {\bf R}^{n-1}$. This discrete map, along with the crossing times, contains all information in the full continuous-time dynamics. 

The structure of the dynamics of an $n$-dimensional network may be represented by a state transition diagram, a directed graph on an $n$-cube, where nodes (vertices, labelled by the Boolean vectors, $\tilde{\bf y}$) represent orthants of phase space and edges represent transitions across orthant boundaries. Under Conditions 1 and 2, the flow across orthant boundaries is unambiguous and determines a direction on the edge between the corresponding nodes on the $n$-cube. If a periodic orbit exists for the network, then it must follow a directed cycle of edges on the $n$-cube. The converse is not necessarily true, as we will see below. 

We now list without proof key properties of the cycle map $M$ (Equation~\ref{eq:Acycle}) and corresponding periodic orbits.

\begin{proposition} Trajectories starting at different points on a given ray through the origin remain on a common ray under iteration of $M$ (though the time taken on the two continuous trajectories will differ). Under Condition 2, such trajectories converge as $t\rightarrow\infty$ (though this convergence may be to the origin).
\end{proposition}
\begin{proposition} Linear subspaces are mapped to linear subspaces by $M$. In particular, straight lines are mapped to straight lines, and planes are mapped to planes.
\end{proposition}

Along a cycle on the $n$-cube, there may be branching nodes, {\em i.e.}, nodes with more than one outgoing edge. These correspond to orthants from which trajectories can exit by more than one boundary hyperplane, depending on which variable reaches zero first. Alternate exit variables impose constraints on the region of an orthant boundary that maps forwards through a specified sequence of boundaries. These constraints take the form of linear inequalities, and the restricted regions are the interiors of `proper cones' (Berman \& Plemmons [\ref{bp94}], p.6). 

\begin{proposition}
Given an $n$-cube cycle and initial orthant boundary, ${\cal O}$, the cone from which trajectories follow the cycle and return to ${\cal O}$ is given by $C=\{{\bf y}\in {\cal O}|R{\bf y}\ge 0\}$, where $R$ is a matrix with one row for each alternate exit variable, $y_i^{(k)}$, around the cycle, each row being 
\beq R_{i,\cdot}=-\frac{{\bf e}_i^{\prime}}{f_i^{(k)}}B^{(k)}B^{(k-1)}\ldots B^{(0)}\,.\label{eq:retconerow}\eeq
\end{proposition}
We allow equality, $R{\bf y}=0$, (trajectories for which two variables cross simultaneously) as limiting cases. Many of the inequalities generated by Equation~\ref{eq:retconerow} will be redundant and can be weeded out in computation.

The domain of definition of $M$ is only $C\subset {\cal O}$. Trajectories starting outside of $C$, but in ${\cal O}$, eventually branch away from the given cycle. Note also that $M$ maps $C$ into ${\cal O}$, not necessarily into $C$. However, a fixed point of the map lying inside $C$ continues to return and corresponds to a periodic orbit for the differential equations. If $C$ is empty, no periodic orbit corresponding to this $n$-cube cycle exists.

\begin{proposition}
Any non-zero (real) fixed point of $M$ (Equation~\ref{eq:Acycle}) in $C$ is a (real) eigenvector of $A$ with eigenvalue $>1$. Conversely, if ${\bf v}$ is a real eigenvector of $A$ with eigenvalue $\lambda >1$, and ${\bf v}\in C$, then 
\beq {\bf y}^*=\frac{(\lambda -1){\bf v}}{\langle\phi ,{\bf v}\rangle }\label{eq:fp}\eeq
is a fixed point of $M$, unique in the span of ${\bf v}$.
If $\lambda =1$, then the only fixed point in the span of ${\bf v}$ is ${\bf 0}$.
\end{proposition}

\begin{proposition} A fixed point, ${\bf y}^*_i$, of $M$ corresponding to the eigenvalue $\lambda_i$ of $A$, is asymptotically stable if $\lambda_i>|\lambda_j|,\,\forall j\ne i$, neutrally stable if $\lambda_i\ge |\lambda_j|,\,\forall j\ne i$, but equality holds for some $j$, and unstable otherwise.
\end{proposition}

\begin{proposition}
A periodic orbit with cycle map $M$ has period $P=\log (\lambda )$, where $\lambda$ is the eigenvalue of the matrix $A$ associated with the fixed point on the orbit.
\end{proposition}

\bigskip
\centerline{\bf 4. APERIODIC BEHAVIOR}
\medskip
\noindent Aperiodic behavior in Glass networks appears common for large $n$. Numerical evidence for ergodicity via an invariant measure for a 6-dimensional network was given by Lewis \& Glass [\ref{lg92}]. For this and a similar 6-dimensional network discussed by Edwards {\em et al} [\ref{ebg99}], aperiodic behavior appeared for certain values of a parameter. Mestl {\em et al} [\ref{mpo95}] showed that chaos cannot exist in 3-dimensional Glass networks, but a 4-dimensional network with a special set of parameters was shown by Mestl {\em et al} [\ref{mlg96}] to have aperiodic dynamics, though it was not proven that there was a chaotic attractor. None of these examples were Boolean Glass networks and it was not known whether 4-dimensional networks of this class could exhibit aperiodicity. 

Numerical experiments that we performed on randomly generated 4-dimensional Boolean Glass networks produced examples for which no fixed point or periodic orbit was detected. We investigate one of these, defined by the interaction function, ${\bf F}$, in Fig.~\ref{fig:chaoscube}. One way to express this network as a system of ODEs is as follows:
\[\begin{array}{lcl}
\dot{y}_1 & = & -y_1+2[\tilde{y}_3]-1 \\
\dot{y}_2 & = & -y_2+2[1-\tilde{y}_3+\tilde{y}_1\tilde{y}_3-\tilde{y}_1\tilde{y}_3\tilde{y}_4]-1 \\
\dot{y}_3 & = & -y_3+2[(1-\tilde{y}_1)(1-\tilde{y}_4)+\tilde{y}_2\tilde{y}_4]-1 \\
\dot{y}_4 & = & -y_4+2[(1-\tilde{y}_1)(1-\tilde{y}_3)+\tilde{y}_1\tilde{y}_2]-1 \,. 
\end{array}\]
A projection of an example 4-dimensional trajectory is shown in Fig.~\ref{fig:chaoscones}a.

\begin{table}[t]
\footnotesize
\begin{tabular}{|c|rrrr|} \hline
orthant ($\tilde{\bf y}$) & \multicolumn{4}{c|}{focal pt. (${\bf F}$)}\\ \hline
0 0 0 0 & $-1$ &  $1$ &  $1$ &  $1$ \\
0 0 0 1 & $-1$ &  $1$ & $-1$ &  $1$ \\
0 0 1 0 &  $1$ & $-1$ &  $1$ & $-1$ \\
0 0 1 1 &  $1$ & $-1$ & $-1$ & $-1$ \\
0 1 0 0 & $-1$ &  $1$ &  $1$ &  $1$ \\
0 1 0 1 & $-1$ &  $1$ &  $1$ &  $1$ \\
0 1 1 0 &  $1$ & $-1$ &  $1$ & $-1$ \\
0 1 1 1 &  $1$ & $-1$ &  $1$ & $-1$ \\ 
1 0 0 0 & $-1$ &  $1$ & $-1$ & $-1$ \\
1 0 0 1 & $-1$ &  $1$ & $-1$ & $-1$ \\
1 0 1 0 &  $1$ &  $1$ & $-1$ & $-1$ \\
1 0 1 1 &  $1$ & $-1$ & $-1$ & $-1$ \\
1 1 0 0 & $-1$ &  $1$ & $-1$ &  $1$ \\
1 1 0 1 & $-1$ &  $1$ &  $1$ &  $1$ \\
1 1 1 0 &  $1$ &  $1$ & $-1$ &  $1$ \\
1 1 1 1 &  $1$ & $-1$ &  $1$ &  $1$ \\ \hline
\end{tabular}
\label{tab:chaoscube}
\end{table}

\begin{figure}[t]
 \vspace{0.15in}
 \special{wmf:c:/document/papers/glassnet/dsa99/chaoscub.wmf x=6.4in y=4.5in}
\caption{\small 4-cube structure of a Boolean Glass network with chaotic behavior. The cycle marked with bold lines corresponds to an unstable periodic orbit.
\label{fig:chaoscube}} 
\end{figure}

\begin{figure}[t]
 \vspace{2.5in}
 \special{wmf:c:/document/papers/glassnet/dsa99/chaos2.wmf x=6in y=4in}
 \caption{\small (a) Projection onto the $y_2$--$y_4$ plane of a trajectory for the 4-dimensional network of Fig.~\ref{fig:chaoscube}. The last 500 of 1000 orthant boundary transitions are shown. (b) Projections of returning cones for two cycles and their images. The triangle $(-1,0),(0,1),(0,0)$ is a projection of the orthant boundary $(0+-+)$. The regions indicated by dotted lines, labelled $C_0$ and $C_1$, are the returning cones.  The regions indicated by solid lines cross-cutting $C_0$ and $C_1$ are the images of these cones under one iteration of their respective maps, $M_0$ and $M_1$. $M_1(C_1)$ (unlabelled) is the narrower region, containing two marked points. The diamonds represent eigenvectors of $M_0$ and $M_1$ and the crosses represent eigenvectors of the composite maps, $M_0M_1$ and $M_1M_0$.
  \label{fig:chaoscones}}
\end{figure}

Consider the two cycles,
\[\begin{array}{l} 
0101\rightarrow 0111\rightarrow 1111\rightarrow 1011\rightarrow {\bf 1001}\rightarrow 1000\rightarrow 1100\rightarrow 1101\quad\mbox{ and} \\
0101\rightarrow 0111\rightarrow 1111\rightarrow 1011\rightarrow {\bf 1010}\rightarrow 1000\rightarrow 1100\rightarrow 1101\,.
\end{array}\]

\noindent Both are feasible and starting from the $(0+-+)$ boundary ({\em i.e.}, the boundary between 1101 and 0101) they have return maps $M_0$ and $M_1$ defined respectively by
\[A_0=\left( \begin{array}{rrr}
     1 & 0 & 0\\
     -2 &  5 & 2\\
     0 & 2 & 1
\end{array}\right)\,,\quad \phi_0=(4,-4,0)^T\,,\]
\[A_1=\left( \begin{array}{rrr}
     1 & -2 & -2\\
     -2 &  -3 & -6\\
     0 & -2 & -3
\end{array}\right)\,,\quad \phi_1=(4,-4,0)^T\,.\]
Their eigenvalues ($\lambda_i$) and corresponding eigenvectors (${\bf v}_i$) are for $A_0$, 
\[\lambda_1\approx 5.8284\,,\quad\lambda_2=1.0000\,,\quad\lambda_3\approx 0.1716\]
\[{\bf v}_1\approx \left( \begin{array}{r}
     0.0000\\
     0.7071\\
     0.2929
\end{array}\right)\,,\quad 
{\bf v}_2=\left( \begin{array}{r}
     0.5000\\
     0.0000\\
     0.5000
\end{array}\right)\,,\quad 
{\bf v}_3\approx \left( \begin{array}{r}
     0.0000\\
     -0.2929\\
     0.7071
\end{array}\right)\,,\]
and for $A_1$,
\[\lambda_1\approx -6.8709\,,\quad\lambda_2\approx 1.9457\,,\quad\lambda_3\approx -0.0748\]
\[{\bf v}_1\approx \left( \begin{array}{r}
     0.2026\\
     0.5257\\
     0.2716
\end{array}\right)\,,\quad 
{\bf v}_2\approx \left( \begin{array}{r}
     0.4728\\
     -0.3754\\
     0.1518
\end{array}\right)\,,\quad 
{\bf v}_3\approx \left( \begin{array}{r}
     -0.2590\\
     -0.4401\\
     0.3009
\end{array}\right)\,.\]
Neither $A_0$ nor $A_1$ has its dominant eigenvector in $(+-+)$, and therefore, neither cycle has a stable periodic orbit. Fig.~\ref{fig:chaoscones}b shows the returning cones for these two cycles in the $(0+-+)$ boundary. In order to depict these 3-dimensional objects in a plane figure, we have projected the cones onto the plane $y_2-y_3+y_4=1$, which is the part of the unit ${\it l}^1$ ball in ${\bf R}^3$ that lies in the $(+-+)$ octant, and then plotted $y_3$ vs. $y_4$. Trajectories starting in the region labelled $C_0$ follow the first cycle above, $M_0$, and return to the $(0+-+)$ boundary. Similarly, trajectories starting in $C_1$ follow $M_1$. 


The images of these two regions under their respective maps are also shown in Fig.~\ref{fig:chaoscones}b. The marked point at $\left( 0,\frac{1}{2}\right)\in C_0$ represents the eigenvector ${\bf v}_2$ of $A_0$ ($\lambda_2>0$) so the ray through ${\bf v}_2$ in ${\bf R}^3$ is invariant under $M_0$ (and thus so is its projection in the figure), but since the corresponding eigenvalue is 1, no non-zero point on the ray is actually fixed and the origin attracts. The stretching and contraction in the directions of the other two eigenvectors (which are perpendicular and lie in the $y_3$--$y_4$ plane) are clearly visible. Since both eigenvalues are positive, there is no inversion of the image in either direction. Note that the eigenvector ${\bf v}_3$ with projected coordinates $(\frac{\sqrt{2}}{2}-1,\frac{\sqrt{2}}{2})\approx (-0.2929,0.7071)$ lies outside $C_0$.

The returning region $C_1$ for map $M_1$ is also subject to stretching and contracting in similar directions. The eigenvector ${\bf v}_2$ associated with eigenvalue $\lambda_2\approx 1.9457$ lies inside $C_1$, so the ray through this point is invariant, and there is an unstable fixed point of $M_1$ on this ray at $(0.1318,-0.1046,0.0423)$, which appears on the projection at $(-0.3754,0.1518)$. The period of the corresponding unstable orbit is $\log (\lambda_2)\approx 0.6656$. The other two eigenvalues are negative, so in the other two directions (in which the $C_1$ region is stretched and contracted) we also have inversion. 

The combined map, defined by $M_0$ and $M_1$, restricted to the projected plane, therefore contains something similar to a Smale horseshoe [\ref{d89}]. It is not exactly topologically equivalent to the horseshoe, but retains many of its properties. Points getting mapped out of $C_0\cup C_1$ will go elsewhere but there is nevertheless a Cantor set, $\Lambda$, in the projected plane, consisting of points whose trajectories remain in $C_0\cup C_1$ both forwards and backwards in time, and an infinite set of unstable periodic points. The main difference from the Smale horseshoe is that no points in the region $M_1(C_1)\cap C_0$ are mapped (by $M_0$) into $C_0$ (this is easy to check by finding the images of the vertices of $M_1(C_1)\cap C_0$ and joining them up by straight lines to find the images of the boundary edges, since straight lines are mapped to straight lines, even when projected). Thus, no points of $\Lambda$, aside from $\left(\frac{1}{2}, 0,\frac{1}{2}\right)$ itself, lie in $M_0(C_0)\cap C_0$ and trajectories of periodic points never follow the $M_0$ cycle twice in a row, {\em i.e.}, their symbolic trajectories do not contain the string `00'. Compositions of $M_0$ and $M_1$ for which $M_0$ is not repeated twice produce unstable periodic orbits. For example, the fixed points corresponding to $M_0M_1$ and $M_1M_0$ are marked in Fig.~\ref{fig:chaoscones}b.

In order to deal with the radial direction that is suppressed in the above projection, we need only confirm that trajectories from points in the Cantor set, $\Lambda$, (or on the rays through points of $\Lambda$) do not converge to the origin. In this case convergence in the radial direction ensures that asymptotically we can ignore the radial component and we approach the 2-dimensional chaotic dynamics of $\Lambda$. 

Boundedness away from the origin can be shown by finding a neighbourhood of the origin in which points on rays through $\Lambda$ always move away from the origin under the map. Simple but tedious calculations show that for each of the corners, $Q_i,\, i=1,\ldots ,12$, of the three regions $M_0(C_0)\cap C_1$, $M_1(C_1)\cap C_1$ and $M_1(C_1)\cap C_0$, the ${\it l}^1$ norm of the appropriate map, $M_j$, satisfies $\| M_j(kQ_i)\|_1>\| kQ_i\|_1=k$ if $0<k<\frac{3}{22}$. Since planes are mapped to planes by each of $M_0$ and $M_1$, the same is true for rays through any point in the interior of one of these three regions (which contain $\Lambda$). Thus, points closer to the origin than $\| {\bf y}\|_1 = \frac{3}{22}$ move away and there is no possibility of convergence to the origin on any ray through $\Lambda$.

\bigskip
\centerline{\bf 5. DISCUSSION}
\medskip
\noindent We have shown that chaotic dynamics exist for the network of Fig.~\ref{fig:chaoscube}. We have not shown that orbits in $\Lambda$ are attracting and, in fact, numerical evidence suggests that they are not. When these equations are integrated, trajectories {\em do} return repeatedly to the orthant boundary $(0+-+)$, but have itineraries including several other cycles besides $M_0$ and $M_1$. These other cycles have returning cones in the gaps left by $C_0$ and $C_1$. Thus, the dynamics are actually more complicated than suggested by the analysis above.

It is surprising that complex dynamics are possible in such a simple network of only 4 variables with Boolean interactions. It certainly is not possible in 4-element discrete-time switching networks. Nevertheless, the techniques discussed here allow considerable progress in analysis of these Glass networks. An important unsolved question is to find necessary or sufficient conditions on the connection structure for chaotic dynamics to occur.

\bigskip
\centerline{\bf ACKNOWLEDGEMENTS}
\medskip
\noindent This work was partially supported by grants from the University of Victoria and the Natural Sciences and Engineering Research Council of Canada.

\bigskip
\centerline{\bf REFERENCES}
\medskip
\small

\newcounter{refno}
\begin{list}
{\arabic{refno}.}{\usecounter{refno}\setlength{\leftmargin}{0.7cm}}
\item\label{bp94} Berman, A., \& Plemmons, R. J., Nonnegative Matrices in the Mathematical Sciences. Academic Press, New York, 1994.
\item\label{d89} Devaney, R. L., An Introduction to Chaotic Dynamical Systems. Addison-Wesley, New York, 1989.
\item\label{ebg99} Edwards, R., Beuter, A. \& Glass, L. Parkinsonian tremor and simplification in network dynamics. Bull. Math. Biol., Vol. 61 (1999) pp. 157--177.
\item\label{g75a} Glass, L., Classification of biological networks by their qualitative dynamics. J. Theor. Biol., Vol. 54 (1975a) pp. 85--107.
\item\label{g75b} Glass, L., Combinatorial and topological methods in nonlinear chemical kinetics. J. Chem. Phys., Vol. 63 (1975b) pp. 1325--1335.
\item\label{g77} Glass, L., Combinatorial aspects of dynamics in biological systems, in: Landman, U. (ed.), Statistical Mechanics and Statistical Methods in Theory and Application. Plenum, New York, 1977, pp.585--611.
\item\label{gh98} Glass, L., \& Hill, C., Ordered and disordered dynamics in random networks. Europhys. Lett., Vol. 41 (1998) pp. 599--604.
\item\label{gp78a} Glass, L., \& Pasternack, J. S.,  Prediction of limit cycles in mathematical models of biological oscillations. Bull. Math. Biol., Vol. 40 (1978a) pp. 27--44.
\item\label{gp78b} Glass, L., \& Pasternack, J. S.,  Stable oscillations in mathematical models of biological control systems. J. Math. Biol., Vol. 6 (1978b) pp. 207--223.
\item\label{lg92} Lewis, J. E., \& Glass, L.,  Nonlinear dynamics and symbolic dynamics of neural networks. Neural Computation, Vol. 4 (1992) pp. 621--642.
\item\label{mbg97} Mestl, T., Bagley, R. J., \& Glass, L.,  Common chaos in arbitrarily complex feedback networks. Phys. Rev. Lett., Vol. 79 (1997) pp. 653--656.
\item\label{mpo95} Mestl, T., Plahte, E., \& Omholt, S.W.,  Periodic solutions in systems of piecewise-linear differential equations. Dynamics and Stability of Systems, Vol. 10 (1995) pp. 179--193.
\item\label{mlg96} Mestl, T., Lemay, C., \& Glass, L.,  Chaos in high-dimensional neural and gene networks. Physica D, Vol. 98 (1996) pp. 33--52.
\end{list}

\vfill\eject
\end{document}